\title{~~\\ On the average number of elements in a finite field with
order or index in a prescribed residue class}
\author{Pieter Moree}
\documentclass[12pt]{article}
\usepackage{amssymb, latexsym, amsfonts}
\def\@ptsize{2}
\newtheorem{Thm}{Theorem}

\newtheorem{Lem}{Lemma}
\newtheorem{Cor}{Corollary}

\newtheorem{Prop}{Proposition}
\newcommand{\qed}{\hfill $\Box$}

\begin{document}
\date{}
\maketitle
\begin{abstract}
\noindent For any prime $p$ the density of elements in
$\mathbb F_p^*$ having
order, respectively index, congruent to $a({\rm mod~}d)$ is being considered.
These densities on average are determined, where
the
average is taken over all finite fields of prime order. Some connections between
the two densties are established. It is also
shown how to compute these densities with high numerical accuracy.\\

\end{abstract}
\section{Introduction}
Let $\mathbb F_q^*$ be the multiplicative group of a finite field
$\mathbb F_q$ and let $x\in \mathbb F_q^*$. The order of $x$ is
the smallest positive integer $k$ such that $x^k=1$, the index is the
largest number $t$ such that $x^{(q-1)/t}=1$. Note that $t=[\mathbb F_q^*:\langle x\rangle]$. Let $\pi(x)$ denote the
number of primes $p\le x$. 
Let $\delta(p;a,d)$ and $\rho(p;a,d)$ denote the density of elements in
$\mathbb F_p^*$ having order, respectively index, congruent to
$a({\rm mod~}d)$.
It is not so difficult to show that
both $\lim_{x\rightarrow \infty}N(a,d)(x)/\pi(x)$ and $\lim_{x\rightarrow \infty}
N'(a,d)(x)/\pi(x)$ exist, where
$$N(a,d)(x)=\sum_{p\le x}\delta(p;a,d){\rm ~and~}
N'(a,d)(x)=\sum_{p\le x}\rho(p;a,d)$$
These limits are denoted by $\delta(a,d)$, respectively
$\rho(a,d)$. In this note these quantities are investigated. The
following result is characteristic of the type of results that will be established. (If 
$a$ and $b$ are natural numbers, then by $(a,b)$, respectively
$[a,b]$, the greatest common divisor, respectively the lowest common multiple of $a$ and $b$ are denoted.
By $\gamma(a)=\prod_{p|a}p$ the square free kernel of $a$ is denoted. Throughout the letter $p$ will
be used to indicate primes.)
\begin{Thm}
\label{drie} $~$\\
{\rm 1)} For every $B>0$ one has
$$N(a,d)(x)=\delta(a,d){\rm Li}(x)+O_B\left({x\over \log^B x}\right),$$
where
$$
\delta(a,d)=\sum_{r=1\atop (1+ra,d)=1}^{\infty}
\sum_{m=1\atop (m,d)|a}^{\infty}{\mu(m)\over mr\varphi(r[m,d])},
$$
and the implied constant depends at most on $B$.\\
{\rm 2)} One has
$$
\delta(a,d)={\varphi((a,d))\over (a,d)}
\sum_{r=1\atop (1+ra,d)=1}^{\infty}\sum_{m=1\atop (m,d)=1}
^{\infty}{\mu(m)\over mr\varphi(mrd)}.
$$
{\rm 3)} If $d|d_1$ and $\gamma(d_1)=\gamma(d)$, then
$\delta(a,d_1)={d\over d_1}\delta(a,d)$.\\
{\rm 4)} One has 
$$\delta(0,d)={1\over d}\prod_{p|d}{1\over 1-{1\over p^2}}{\rm ~and~}
\delta(d,2d)=\cases{{1\over 2}\delta(0,2d) &if $d$ is odd;\cr
\delta(0,2d) &if $d$ is even.}$$
{\rm 5)} For $s\ge 1$ one has
$$\delta(a,2^s)=\cases{2^{2-s}/3 &if $a$ is even;\cr
2^{1-s}/3 &if $a$ is odd.}$$
{\rm 6)} Let $q$ be a prime and $q\nmid a$. Then
$$\delta(a,q)={q^2\over (q-1)(q^2-1)}-{q\over q^2-q-1}\rho(-{1\over a},q).$$
{\rm 7)} Put $W_d(a)=\{0\le r<d~:~(1+ra,d)=1\}$.
One has
$$\delta(a,d)={\varphi((a,d))\over (a,d)\varphi(d)\prod_{p|d}(1-{1\over p(p-1)})}
\sum_{\alpha\in W_d(a)}\rho(\alpha,d)\prod_{p|(\alpha,d)}{p^2-p-1\over p^2-1}.$$
\end{Thm} 
In the next subsection a characteristic zero version of $\delta(a,d)$ and $\rho(a,d)$ will
be discussed. Indeed, these characteristic zero quantities (exhibiting far
more complicated behaviour) motivated the author
to study $\delta(a,d)$ and $\rho(a,d)$. The behaviour of these characteristic zero
quantities turns out to have many resemblances with that of $\delta(a,d)$ and $\rho(a,d)$.
\subsection{Connections with characteristic zero}
Let $g\in \mathbb Q\backslash \{-1,0,1\}$ and $p$ be a prime. By $\nu_p(g)$ the exponent of
$p$ in the canonical factorisation of $g$ is denoted. If $\nu_p(g)=0$, then $g$
can be considered as an element of $\mathbb F_p^*$ with order ord$_g(p)$ and index $r_g(p)$. 
Let $N_g(a,d)$ and $N'_g(a,d)$ denote the set of primes $p$ with
$\nu_p(g)=0$
such that the order, respectively index of $g({\rm mod~}p)$ is congruent to
$a({\rm mod~}d)$.\\ 
\indent In 
case $g=2$ the set
$N'_g(a,d)$ was first considered by Pappalardi \cite{Pappalardi}, 
who proved that it has
a natural density $\rho_g(a,d)$ under 
the assumption of the Generalized Riemann Hypothesis (GRH). For
the general case see \cite{Moreealleen}.\\
\indent The methods of \cite{Moreealleen} can be extended (see 
\cite{Moreealleen2}) to show that under
GRH the set
$N_g(a,d)$ has a natural density $\delta_g(a,d)$,
the evaluation of which seems to be far less easy 
than that of $\rho_g(a,d)$. By ${\bar \delta}_g(d)$ the vector
$(\delta_g(0,d),\cdots,\delta_g(d-1,d))$ (if it exists) is denoted. Up to
this century only ${\bar \delta}_g(2)$ had been evaluated.
Recently Chinen and Murata \cite{CM4,CM5} computed
${\bar \delta}_g(4)$ (on GRH) under the assumption that $g$ is a positive
integer that is not a pure power. In \cite{Moreealleen}, on GRH,
${\bar \delta}_g(3)$ and ${\bar \delta}_g(4)$ are evaluated for
each $g\in \mathbb Q\backslash\{-1,0,1\}$.
If $d|2(a,d)$, then
$\delta_g(a,d)$ can be
evaluated unconditionally, cf. \cite{Ballot, Odoni, Wiertelak1, Wiertelak2}.\\
\indent Let $G$ be the set of rational numbers $g$ that cannot be written
as $-g_0^h$ or $g_0^h$ with $h>1$ an integer and $g_0$ a rational
number. By $D(g)$ the discriminant of the number field
$\mathbb Q(\sqrt{|g|})$ is denoted.
The functions $N(a,d)(x)$ and $N'(a,d)(x)$ can be considered as 
naive heuristic approximations of $N_g(a,d)(x)$ and $N'_g(a,d)(x)$ (if $S$ is any set of non-negative 
integers then $S(x)$ denotes the number of elements in $S$
not exceeding $x$).
For more on heuristics and primitive root theory, see
\cite{Early, Near}. Theorem \ref{conjeen} and Theorem \ref{een} say that
as $D(g)$ becomes large, the naive heuristic for $\delta_g(a,d)$ and
$\rho_g(a,d)$ become more and more accurate. Some related numerical
material is
provided in Table 3 and Table 4. The next result is proved in \cite{Moreealleen2}.
\begin{Thm} {\rm (GRH)}.
\label{conjeen}
If $D(g)\rightarrow \infty$ with $g\in G$, then
$\delta_g(a,d)$ tends to $\delta(a,d)$.
\end{Thm}
In this paper the following similar (but easier) result will be proved.
\begin{Thm} {\rm (GRH)}.
\label{een}
If $D(g)\rightarrow \infty$ with $g\in G$, then
$\rho_g(a,d)$ tends to $\rho(a,d)$.
\end{Thm}
\indent It turns out that both $\delta(a,d)$ and $\rho(a,d)$ are much more
accessible quantities than $\delta_g(a,d)$, respectively $\rho_g(a,d)$.
In the light of
the latter two theorems it 
thus seems of some importance to compute $\delta(a,d)$
and $\rho(a,d)$, which is the purpose of this paper. The more
complicated nature of $N_g(a,d)$ versus $N_g'(a,d)$ is mirrored in
the fact that $\delta(a,d)$ is rather more difficult to compute than
$\rho(a,d)$.\\
\indent As a prelude to proving Theorem \ref{een}, its
`index equals $t$' analog is proved in Section \ref{indext}.  
Theorem \ref{een} is then proved in Section 3 (this involves
evaluating $\rho(a,d)$).
In Section 4 Theorem \ref{drie} is considered.
In Section 5 an Euler product $A_{\chi}$ involving a Dirichlet character $\chi$
is studied and it is shown how $\rho(a,d)$ and $\delta(a,d)$ can be
expressed in terms of $A_{\chi}$'s. Since $A_{\chi}$
can be evaluated with high numerical accuracy (Sections 6 and 7) this then allows
us to evaluate $\rho(a,d)$ and $\delta(a,d)$
with high numerical precision.

\section{Index $t$}
\label{indext}
Let $\rho_p(t)$ denote the density of elements in $\mathbb F_p^*$ having
index $t$. Note that $\rho_p(t)=\varphi((p-1)/t)/(p-1)$ if $p\equiv
1({\rm mod~}t)$ and $\rho_p(t)=0$ otherwise. By the method of
proof of Theorem \ref{pappathm} 
(cf. \cite[p. 161]{Near}) it is easy to show that, for every $B>0$, 
\begin{equation}
\label{eerstesom}
\sum_{p\le x}\rho_p(t)=\sum_{p\le x\atop p\equiv 1({\rm mod~}t)}
{\varphi\left({p-1\over t}\right)\over p-1}=
{\rm Li}(x)\sum_{n=1}^{\infty}{\mu(n)\over tn\varphi(tn)}
+O\left({x\over \log^B x}\right),
\end{equation}
where Li$(x)$
denotes the logarithmic integral.
Now
$$\sum_{n=1}^{\infty}{\mu(n)\over tn\varphi(tn)}=
{1\over t\varphi(t)}\sum_{n=1}^{\infty}{\mu(n)\varphi(t)\over n\varphi(tn)}.$$
The latter sum has as argument a multiplicative function 
in $n$. On applying
Euler's identity, it is then inferred that
\begin{equation}
\label{art}
\sum_{n=1}^{\infty}{\mu(t)\over nt\varphi(nt)}=Ar(t),
\end{equation}
where $$A=\prod_p\left(1-{1\over p(p-1)}\right)
=0.37395581361920228805472805434641641511\cdots$$
denotes Artin's constant and 
$$r(t)={1\over t^2}\prod_{p|t}{p^2-1\over p^2-p-1}={1\over t\varphi(t)}
\prod_{p|t}{1-{1\over p^2}\over 1-{1\over p(p-1)}}.$$
Note that
\begin{equation}
\label{er}
{1\over t\varphi(t)}\le r(t) \le {6\over A\pi^2t\varphi(t)}.
\end{equation}
Combination of (\ref{eerstesom}) and (\ref{art}) yields the following result.
\begin{Prop}
For every $B>0$, one has
$$
\sum_{p\le x}\rho_p(t)=Ar(t)
{\rm Li}(x)
+O\left({x\over \log^B x}\right),
$$
Thus, the density of elements in $\mathbb F_p^*$ with index $t$
equals $Ar(t)$ on average.
\end{Prop}
\indent Let $N_g(t)$ denote the set of primes $p$ with $r_g(p)=t$. If
$g\in \mathbb Q\backslash \{-1,0,1\}$, then it can be shown 
\cite{Lenstra, Wagstaff}, under GRH, that
$N_g(t)$ has density
\begin{equation}
\label{wagyourtail}
\delta(N_g(t))=\sum_{n=1}^{\infty}{\mu(n)\over 
[\mathbb Q(\zeta_{nt},g^{1/nt}):\mathbb Q]}.
\end{equation}
For an explicit evaluation of this density see \cite{Murato, Wagstaff}.
We can now prove the following result.
\begin{Prop} 
\label{proposition}
{\rm (GRH)}. Let $g\in G$.
If $g>0$, set $m=[2,D(g)]$.
If $g<0$, set $m=D(g)/2$ if $D(g)\equiv 4({\rm mod~}8)$
and $m=[4,D(g)]$ otherwise. Put $m_1=m/(t,m)$.
If $g\in G$, then
$$|\delta(N_g(t))
-Ar(t)|\le {2.21\over tm_1\varphi(tm_1)}.$$
\end{Prop}
\begin{Cor}
{\rm (GRH)}.
The density of $N_g(t)$ exists and
if $D(g)\rightarrow \infty$ with $g\in G$, tends
to the average density, $Ar(t)$, of elements in $\mathbb F_p^*$ having
index $t$.
\end{Cor}
The latter corollary is the 
`index equals $t$' analog of Theorem \ref{een}.\\

\noindent {\it Proof of Proposition} \ref{proposition}. By (\ref{wagyourtail})
and the evaluation of the degree $[\mathbb Q(\zeta_k,g^{1/k}):\mathbb Q]$
as given in \cite{Wagstaff}, it is deduced that
\begin{equation}
\label{komiknogopterug}
\delta(N_g(t))=Ar(t)+\sum_{k=1\atop m|kt}^{\infty}{\mu(k)\over kt\varphi(kt)}
=Ar(t)+\sum_{k=1\atop m_1|k}^{\infty}{\mu(k)\over kt\varphi(kt)}.
\end{equation}
On noting that $\varphi(zw)\ge \varphi(w)\varphi(z)$, with $w$ and $z$
arbitrary integers and using that $\sum_{k}1/(k\varphi(k))<2.21$, the
result then follows. \qed\\

\noindent {\tt Remark 1}. The sum $\sum_{k}1/(k\varphi(k))$ can be written as an Euler product
of the form $\prod_p F_1(p)/F_2(p)$, with $F_j(X)\in \mathbb Z[X]$ for $j=1,2$ and monic.
Using Theorem 2 of \cite{Moreeconstant} such Euler products 
can be expressed in terms of values at integer points of the (partial) Riemann zeta-function.
This enables one to evaluate these constants with hunderds of decimals of precision, see 
\cite{Niklasch}. A similar
idea forms the basis of Theorem \ref{numeriek} and Theorem \ref{fiboconv}.

\section{Computation of $\rho(a,d)$}
Equation (\ref{wagyourtail}) suggests that, 
under GRH, one should have
\begin{Prop} {\rm (GRH)}. If $g\in \mathbb Q\backslash \{-1,0,1\}$, then
\label{pappa}
$$
\rho_g(a,d)=\sum_{t\equiv a({\rm mod~}d)}\delta(N_g(t))=\sum_{t\equiv a({\rm mod~}d)}
\sum_{n=1}^{\infty}{\mu(n)\over 
[\mathbb Q(\zeta_{nt},g^{1/nt}):\mathbb Q]}.
$$
\end{Prop}
(In this proposition and in the sequel sums over $t$ are assumed to run over
positive integers only.)
Indeed by \cite{Pappalardi}, cf. \cite{Moreealleen}, Proposition \ref{pappa} is known to be true.
Similarly one would expect that $\rho(a,d)$ satisfies (\ref{rhoaccent}) as can indeed
be proved.
\begin{Thm} 
\label{pappathm}
For every $B>0$ one has
\begin{equation}
\label{pappathmeq}
N'(a,d)(x)=\rho(a,d){\rm Li}(x)+O\left({x\over \log^B x}\right),
\end{equation}
where
\begin{equation}
\label{rhoaccent}
\rho(a,d)=A\sum_{t\equiv a({\rm mod~}d)}r(t)=
\sum_{t\equiv a({\rm mod~}d)}{A\over t^2}\prod_{p|t}{p^2-1\over p^2-p-1}
\end{equation} 
and the implied constant
depends at most on $B$.
\end{Thm}
{\it Proof}. One has
$$N'(a,d)(x)=\sum_{p\le x}\sum_{t|p-1\atop t\equiv a({\rm mod~}d)}
{\varphi({p-1\over t})\over p-1}.$$
On using that $\varphi(n)/n=\sum_{m|n}\mu(m)/m$, one obtains
$$N'(a,d)(x)=\sum_{p\le x}\sum_{t|p-1\atop t\equiv a({\rm mod~}d)}
\sum_{m|{p-1\over t}}{\mu(m)\over mt}.$$
Writing $v=mt$ and bringing the summation over $p$ to the inside, one
obtains
$$N'(a,d)(x)=\sum_{v\le x-1}{1\over v}\sum_{t|v\atop t\equiv a({\rm mod~}d)}\mu({v\over t})
\sum_{p\le x\atop p\equiv 1({\rm mod~}v)}1.$$
The summation range is split up into the range
$v\le \log^{B+1}x$ and $\log^{B+1} x<v\le x$. In the former range
the Siegel-Walfisz theorem (see e.g. \cite[Satz 4.8.3]{Prachar}) is invoked
and for the latter range the trivial estimate
$\sum_{p\le x,~p\equiv 1({\rm mod~}v)}1<x/v$ is employed. 
Let $d(v)$ denote the number of divisors of $v$.
Together with
the trivial estimate
$|\sum_{t|v,~t\equiv a({\rm mod~}d)}\mu(v/t)|\le d(v)\ll v^{\epsilon}$,
which holds for every $\epsilon>0$, one concludes (cf. \cite[p. 161]{Near}) that
(\ref{pappathmeq}) holds with
\begin{equation}
\label{nogonbedorven}
\rho(a,d)=\sum_{v=1}^{\infty}{\sum_{t|v,~t\equiv a({\rm mod~}d)}\mu({v\over t})
\over v\varphi(v)}.
\end{equation}
Interchanging the order of summation and using (\ref{art}) one infers that
$$
\rho(a,d)=\sum_{t\equiv a({\rm mod~}d)}\sum_{v_1=1}^{\infty}{\mu(v_1)\over tv_1\varphi(tv_1)}
=A\sum_{t\equiv a({\rm mod~}d)}r(t).
$$
This concludes the proof. \qed\\

\noindent Let $a>0$. As $d$
becomes large, the first term in the second summation in
(\ref{rhoaccent}), $Ar(a)$, tends to be dominant by Corollary \ref{BLAH}.
In particular, $\lim_{d\rightarrow \infty}\rho(a,d)=Ar(a)$.

\begin{Prop}
\label{null}
One has $$\rho(0,d)={1\over d\varphi(d)} {\rm ~and~}
\rho(d,2d)=\cases{\rho(0,2d) &if $d$ is odd;\cr
3\rho(0,2d) &if $d$ is even.}$$
\end{Prop}
{\it Proof}. Note that
$$\rho(0,d)=r(d)\sum_{m=1}^{\infty}{r(dm)\over r(d)}.$$
Since $r(dm)/r(d)$ is a multiplicative function in $m$, the
identity for $\rho(0,d)$ then follows on applying Euler's identity and
noting that $\varphi(d)/d=\prod_{p|d}(1-1/p)$. 
The identity for $\rho(0,d)$ together with the observation that
$\rho(0,2d)+\rho(d,2d)=\rho(0,d)$, then yields the
truth of the remainder of the assertion. \qed\\

\noindent By $\omega(m)$ the number of distinct prime
divisors of $m$ is denoted.
\begin{Prop} {\rm (GRH)}.
\label{mamma}
Let $g\in G$ and $m$ be as in Proposition {\rm \ref{proposition}}, then
$$\Big|\rho_g(a,d)-\rho(a,d)\Big|\le {2^{\omega(m)+2}\over m\varphi(m)}.$$
\end{Prop}
{\it Proof}. By Theorem \ref{pappa}, 
Proposition \ref{pappa} and (\ref{komiknogopterug}) one infers 
on putting $kt=v$ that,
under GRH,
$$\Big|\rho_g(a,d)-\rho(a,d)\Big|\le \sum_{t\equiv a({\rm mod~}d)}
\sum_{k=1\atop m|kt}^{\infty}{|\mu(k)|\over kt\varphi(kt)}
= \sum_{m|v}{\sum_{t|v,~t\equiv a({\rm mod~}d)}
|\mu({v\over t})|\over v\varphi(v)}.$$
On noting that  
$$\sum_{m|v}{\sum_{t|v,~t\equiv a({\rm mod~}d)}
|\mu({v\over d})|\over v\varphi(v)}\le \sum_{v=1}^{\infty}{2^{\omega(mv)}\over mv\varphi(mv)}
\le {2^{\omega(m)}\over m\varphi(m)}\sum_{v=1}^{\infty}{2^{\omega(v)}\over v\varphi(v)}
\le {2^{\omega(m)+2}\over m\varphi(m)},$$
the result follows. \qed\\

\noindent Since $2^{\omega(m)}/(m\varphi(m))$ tends to zero with increasing $m$ and $m$ tends to
infinity as $D(g)$ tends to infinity, Theorem \ref{een} is a  consequence of
Proposition \ref{mamma}.\\
\indent The following result is concerned with $\mathbb Q$-linear relations between
the $\rho(a,d)$'s with $d$ fixed.
\begin{Lem}
\label{vectorspace}
Let $\alpha_1,\cdots,\alpha_{\varphi(d)}$ be representatives of the reduced residue
classes mod $d$. Then, for every integer $a$,
$$\rho(a,d)\in \mathbb Q[\rho(\alpha_1,d),\cdots,\rho(\alpha_{\varphi(d)-1},d)].$$
\end{Lem}
\begin{Cor} If $d|d_1$ and $\beta_1,\cdots,\beta_{\varphi(d_1)}$ are representatives of the reduced residue
classes mod $d_1$, then
$$\mathbb Q[\rho(\alpha_1,d),\cdots,\rho(\alpha_{\varphi(d)-1},d)]\subseteq
\mathbb Q[\rho(\beta_1,d_1),\cdots,\rho(\beta_{\varphi(d_1)-1},d_1)].$$
\end{Cor}
{\it Proof of Lemma} \ref{vectorspace}. It is easy to show that
\begin{equation}
\label{reducedtelop}
\sum_{j=1}^{\varphi(d)}\rho(\alpha_j,d)=A\sum_{(t,d)=1}r(t)=\prod_{p|d}\left(1-{1\over p(p-1)}\right)
\in \mathbb Q.
\end{equation}
It is thus enough to show that
$\rho(a,d)\in V_d:=\mathbb Q[\rho(\alpha_1,d),\cdots,\rho(\alpha_{\varphi(d)},d)]$.\\
\indent Let $\alpha=a/(a,d)$ and $\delta=d/(a,d)$. Note that $(\alpha,\delta)=1$. Let $\delta_1$
be the largest divisor of $(a,d)$ with $(\delta,\delta_1)=1$ and write $(a,d)=\delta_1\delta_2$.
If $\alpha,\delta,\delta_1$ and $\delta_2$ are being used for integers 
other then $a$ and $d$, then this will be made explicit in the 
notation. Thus the meaning of $\delta_1(a_j,d)$, which appears
later in the proof, should be obvious.
Note that
$$\rho(a,d)=A\sum_{t\equiv \alpha({\rm mod~}\delta)}r(\delta_1\delta_2t)=
Ar(\delta_2)\sum_{t\equiv \alpha({\rm mod~}\delta)}r(\delta_1t).$$
The proof proceeds with induction with respect to the
number of distinct prime divisors of $\delta_1$. If $\omega(\delta_1)=0$, then $\delta_1=1$
and one has to show that $\rho(\alpha,\delta)\in V_d$, where $\alpha({\rm mod~}\delta)$ is
a reduced residue class mod $\delta$. Since $\alpha({\rm mod~}\delta)$ lifts to $d/\delta$
reduced residue classes mod $d$, this is clear. If $\omega(\delta_1)=1$, then $\delta_1=q^e$
with $q$ a prime and $e\ge 1$. Then one has
\begin{eqnarray}
\rho(a,d)&=&Ar(\delta_2)\sum_{t\equiv \alpha({\rm mod~}\delta)}r(q^et)\nonumber\cr
&=&r(q^e)Ar(\delta_2)\sum_{t\equiv \alpha({\rm mod~}\delta)\atop t\not\equiv 0({\rm mod~}q)}r(t)
+{A\over q^2}r(\delta_2)\sum_{t\equiv {\alpha\over q}({\rm mod~}\delta)}r(q^et)\nonumber\cr
&=& c_1+{A\over q^2}r(\delta_2)\sum_{t\equiv {\alpha\over q}({\rm mod~}\delta)}r(q^et)\nonumber\cr
&~& ~~: \nonumber\cr
&=& c_n+{A\over q^{2n}}r(\delta_2)\sum_{t\equiv {\alpha\over q^n}({\rm mod~}\delta)}r(q^et),
\end{eqnarray}
where $c_i\in V_d$. On choosing $n\ge 1$ to be such that $q^n\equiv 1({\rm mod~}\delta)$, one infers
that $\rho(a,d)\in V_d$. Suppose the result has been proved for all $a$ and $d$ with
$\omega(\delta_1)\le m$ for some $m\ge 1$. Then consider next $a$ and $d$ with $\omega(\delta_1)=m+1$.
 One has
\begin{equation}
\label{technisch}
\rho(a,d)=\sum_{j=1}^{\delta_1}Ar(\delta_2)\sum_{t\equiv \alpha({\rm mod~}\delta)\atop
t\equiv j({\rm mod~}\delta_1)}r(\delta_1t).
\end{equation}
Note that $$A\sum_{t\equiv \alpha({\rm mod~}\delta)\atop
t\equiv j({\rm mod~}\delta_1)}r(\delta_1t)=A\sum_{t\equiv \alpha({\rm mod~}\delta)\atop
t\equiv {j\over (j,\delta_1)} ({\rm mod~}{\delta_1\over (j,\delta_1)})}r(\delta_2)r(\delta_1(j,\delta_1)t).$$
The latter sum equals a rational multiple times
$$A\sum_{t\equiv \alpha({\rm mod~}\delta)\atop
t\equiv {j\over (j,\delta_1)} ({\rm mod~}{\delta_1\over (j,\delta_1)})}r(\delta_2)r((j,\delta_1)t)
=\rho(a_j,d),$$
for some integer $a_j$. Note that $\delta_1(a_j,d)=(j,\delta_1)$. 
Thus by the induction hypothesis
all terms in (\ref{technisch}) with $\gamma(\delta_1)\nmid j$ are in $V_d$ (where
$\gamma(\delta_1)$ denotes the squarefree kernel of $\delta_1$).
One infers that
\begin{eqnarray}
\rho(a,d)
&=& d_1 + A r(\delta_2)\sum_{t\equiv \alpha({\rm mod~}\delta )\atop t\equiv 0({\rm mod~}\gamma(\delta_1))}r(\delta_1t)\nonumber\cr
&=& d_1+{A\over \gamma(\delta_1)^2}r(\delta_2)\sum_{t\equiv {\alpha\over \gamma(\delta_1)}({\rm mod~}\delta)}
r(\delta_1t)\nonumber\cr
&~& ~~: \nonumber\cr
&=& d_n+{A\over \gamma(\delta_1)^{2n}}r(\delta_2)\sum_{t\equiv {\alpha\over 
\gamma(\delta_1)^n}({\rm mod~}\delta)}r(\delta_1t),
\end{eqnarray}
where $\delta_i\in V_d$. On choosing $n\ge 1$ to be such that $\gamma(\delta_1)^n\equiv 
1({\rm mod~}\delta)$ one infers that $\rho(a,d)\in V_d$. \qed\\

\noindent {\tt Example 1}. The result says that $\rho(a,6)\in \mathbb Q[\rho(1,6)]$. Indeed,
$\rho(0,6)=\rho(3,6)=1/12$. Furthermore $\rho(2,6)=1/12+3\rho(1,6)/5$, $\rho(4,6)=1/3-3\rho(1,6)/5$
and $\rho(5,6)=5/12-\rho(1,6)$.

\section{Proof of Theorem \ref{drie}}
\noindent {\it Proof of Theorem} \ref{drie}. 1) Note that
$$N(a,d)(x)=\sum_{p\le x}\delta(p;a,d)=\sum_{p\le x}\sum_{r|p-1\atop {p-1\over r}\equiv a({\rm mod~}d)}
{\varphi\left({p-1\over r}\right)\over p-1}.$$
Proceeding as in the proof of Theorem \ref{pappathm}, one infers that
\begin{equation}
\label{triplesumm}
N(a,d)(x)=\sum_{r\le x-1}\sum_{m\le {x-1\over r}}{\mu(m)\over mr}
\sum_{{p\le x,~p\equiv 1({\rm mod~}rm)}\atop p\equiv 1+ra({\rm mod~}rd)}1.
\end{equation}
Now for the inner sum to be non-zero the two congruences must be
compatible. By the Chinese remainder theorem this is the case if and
only if $1\equiv 1+ra ({\rm mod~}r(m,d))$, that is if and only if
$a\equiv 0({\rm mod~}(m,d))$. If the two congruences are compatible, then
they form a reduced residue class if and only if $(1+ra,d)=1$. If
the residue class is not reduced it contains at most one prime and
the contribution of these primes to $N(a,d)(x)$ is
bounded in absolute value by $\sum_{v\le x}d(v)/v=O(\log^2 x)$.\\
\indent The summation range $mr\le x-1$ in (\ref{triplesumm}) is split up into the range
$r[m,d]\le \log^C x$ and the range $r[m,d]>\log^C x$, where $C$ is to be chosen
later. All error terms arising in this way are easily seen to be
of the claimed order of growth, except the error term
$$E(x)={\rm Li}(x)\sum_{r}\sum_{(m,d)|a\atop r[m,d]>\log^C x}{|\mu(m)|\over mr\varphi(r[m,d])},$$
which arises on completing the sum
$${\rm Li}(x)\sum_{r\le x-1\atop (1+ra,d)=1}\sum_{m\le {x-1\over r},~(m,d)|a\atop r[m,d]\le \log^C x}
{|\mu(m)|\over mr\varphi(r[m,d])},$$
to $\delta(a,d)$. On noting that $r[m,d]>\log^C x$ implies $rmd>\log^C x$ and using
that $\varphi(zw)\ge \varphi(z)\varphi(w)$, one obtains,
cf. the proof of part 2, 
$$E(x)=O\left( {{\rm Li}(x)\over \varphi(d)}\sum_{m_1|(a,d)}{|\mu(m_1)|\over m_1}
\sum_{r}\sum_{m_2\atop rm_2>\log^C x/(dm_1)}{|\mu(m_2)|\over rm_2\varphi(rm_2)}\right).$$
From this $E(x)$ is easily seen to be $O_B(x/\log^B x)$, when $C$ is chosen to be
sufficiently large.\\
{\rm 2)} By part 1 it is enough to show that
$$I_1:=\sum_{m=1\atop (m,d)|a}^{\infty}{\mu(m)\over mr\varphi(r[m,d])}
={\varphi((a,d))\over (a,d)}\sum_{m=1\atop (m,d)=1}^{\infty}{\mu(m)\over mr\varphi(mrd)}.$$
Note that
$$I_1=\sum_{m_1|(a,d)}\sum_{m=1\atop (m,d)=m_1}^{\infty}{\mu(m)\over mr\varphi({mrd\over m_1})}.$$
On writing $m=m_1m_2$ one obtains
\begin{eqnarray}
I_1&=&\sum_{m_1|(a,d)}\sum_{m_2=1\atop (m_2,d/m_1)=1}^{\infty}{\mu(m_1m_2)\over m_1m_2r\varphi(m_2rd)}\cr
&=&\sum_{m_1|(a,d)}{\mu(m_1)\over m_1}\sum_{m_2=1,~(m_2,d/m_1)=1
\atop (m_2,m_1)=1}^{\infty}
{\mu(m_2)\over rm_2\varphi(m_2rd)}\cr
&=&\sum_{m_1|(a,d)}{\mu(m_1)\over m_1}\sum_{m_2=1\atop (m_2,d)=1}^{\infty}{\mu(m_2)\over rm_2\varphi(m_2rd)}\cr
&=&{\varphi((a,d))\over (a,d)}\sum_{m=1\atop (m,d)=1}^{\infty}{\mu(m)\over rm\varphi(mrd)}.\nonumber
\end{eqnarray}
3) The condition on $d$ and $d_1$ ensures that $(1+ra,d_1)=1$ iff
$(1+ra,d)=1$ and $(m,d_1)=1$ iff $(m,d)=1$. Furthermore one has
$\varphi((a,d_1))/(a,d_1)=\varphi((a,d))/(a,d)$. 
By part 2 one then finds
$$\delta(a,d_1)={\varphi((a,d))\over (a,d)}\sum_{r=1\atop (1+ra,d)=1}^{\infty}
\sum_{m=1\atop (m,d)=1}^{\infty}{\mu(m)\over mr\varphi(mrd_1)}.$$
On noting that $\varphi(mrd_1)=\varphi(mrd)d_1/d$, the proof of part 2 is
completed.\\
4) By part 2 one has, on writing $mr=v$,
$$\delta(0,d)={\varphi(d)\over d}
\sum_{v=1}^{\infty}{\sum_{m|v,~(m,d)=1}\mu(v)\over v\varphi(vd)}.$$
The inner sum equals one if $\gamma(v)|d$ and zero otherwise. Thus
$$\delta(0,d)={\varphi(d)\over d}\sum_{\gamma(v)|d}{1\over v\varphi(vd)}
={1\over d}\sum_{\gamma(v)|d}{1\over v^2}
={1\over d}\prod_{p|d}{1\over 1-{1\over p^2}}.$$
The formula for $\delta(d,2d)$ easily follows from that of
$\delta(0,d)$ and the observation that $\delta(0,2d)+\delta(d,2d)=\delta(0,d)$.\\
5) An easy consequence of part 3 and part 4.\\
6) Using part 2 and
(\ref{art}), one infers that
\begin{eqnarray}
\delta(a,q)&=&\delta(0,q)-\sum_{t\equiv -{1\over a}({\rm mod~}q)}\sum_{m=1\atop (m,q)=1}
^{\infty}{\mu(m)\over mt\phi(mtq)}\nonumber\cr
&=&\delta(0,q)-{1\over q-1} \sum_{t\equiv -{1\over a}({\rm mod~}q)}\sum_{m=1\atop (m,q)=1}
^{\infty}{\mu(m)\over mt\phi(mt)}\nonumber\cr
&=&\delta(0,q)-{1\over q-1} \sum_{t\equiv -{1\over a}({\rm mod~}q)}{A\over 1-{1\over q(q-1)}}r(t).
\end{eqnarray}
The proof is then completed on invoking 
part 3 and (\ref{rhoaccent}). \\
7) By part 2 one can write
\begin{equation}
\label{nogeen}
\delta(a,d)={\varphi((a,d))\over (a,d)}
\sum_{r=1\atop (1+ra,d)=1}^{\infty}{1\over r\varphi(rd)}\sum_{m=1\atop (m,d)=1}
^{\infty}{\mu(m)\varphi(rd)\over m\varphi(mrd)}.
\end{equation}
Denote the inner sum in (\ref{nogeen}) by $I_2$. One has
\begin{eqnarray}
I_2&=&\prod_{p\nmid rd}(1-{1\over p(p-1)})\prod_{p|r\atop p\nmid d}
(1-{1\over p^2})\cr
&=&{A\over \prod_{p|d}(1-{1\over p(p-1)})}\prod_{p|r}
{{1-{1\over p^2}}\over {1-{1\over p(p-1)}}}\prod_{p|(r,d)}{1-{1\over p(p-1)}\over
1-{1\over p^2}}.\nonumber
\end{eqnarray}
One thus obtains that 
$$\delta(a,d)={\varphi((a,d))\over (a,d)\prod_{p|d}(1-{1\over p(p-1)})}\sum_{\alpha\in W_d(a)}
\sum_{w=1\atop w\equiv \alpha({\rm mod~}d)}^{\infty}
{A\varphi(w)\over \varphi(wd)}r(w)\prod_{p|(\alpha,d)}{1-{1\over p(p-1)}\over
1-{1\over p^2}}.$$
Note that
$${\varphi(wd)\over \varphi(w)}=d\prod_{p|d\atop p\nmid w}(1-{1\over p})=
{d\prod_{p|d}(1-{1\over p})\over \prod_{p|(w,d)}(1-{1\over p})}
={\varphi(d)\over \prod_{p|(\alpha,d)}(1-{1\over p})}.$$
This equation together with the latter one derived for $\delta(a,d)$ and (\ref{art}), then
yields the result. \qed

\section{The densities and $A_{\chi}$}
Given a Dirichlet character $\chi$ mod $d$, let $h_{\chi}=\chi\star \mu$,
that is $h_{\chi}$ denotes the Dirichlet convolution of $\chi$ and 
the M\"obius function $\mu$.
Note the following
trivial result.
\begin{Lem}
The function $h_{\chi}$ is multiplicative and satisfies $h_{\chi}(1)=1$
and furthermore $h_{\chi}(p^r)=\chi(p)^{r-1}[\chi(p)-1]$, where the
convention that $0^{0}=1$ is adopted.
\end{Lem}
In particular if $\chi$ is the trivial character mod $d$, then
\begin{equation}
\label{muutje}
h_{\chi}(v)=\cases{\mu(v) &if $v|d$;\cr 0 &otherwise.}
\end{equation}
By using one of the orthogonality relations for Dirichlet characters, the 
following result is easily inferred.
\begin{Lem}
\label{nonmulttomult}
Let $a({\rm mod~}d)$ be a reduced residue class mod $d$. One has
$$\sum_{t\equiv a({\rm mod~}d)\atop t|v}\mu({v\over t})=
{1\over \varphi(d)}\sum_{\chi({\rm mod~}d)}{\overline{\chi(a)}}h_{\chi}(v),$$
where $\chi$ runs over the Dirichlet characters modulo $d$.
\end{Lem}
The reader is referred to Section 2.4 of \cite{Moreealleen} for some 
further properties of $h_{\chi}$.\\
\indent In what follows
sums of the form
$$\sum_{(v,d)=1}{h_{\chi}(v)\over v\varphi(v)}$$
will feature. It is easy to see that this sum is
absolutely convergent. Since its argument is multiplicative, one then
obtains that the latter sum equals
$$A_{\chi}=\prod_{p:\chi(p)\ne
0}\left(1+{[\chi(p)-1]p\over [p^2-\chi(p)](p-1)}\right).$$ 
Note that if $\chi$ is the principal character, then $A_{\chi}=1$.
For a fixed prime $p$ and $\alpha\in \mathbb R$, $0\le \alpha<1$ let
$$f_p(\alpha)=
\left(1+{(e^{2\pi i\alpha}-1)p\over (p^2-e^{2\pi i\alpha})(p-1)}\right).$$
A tedious analysis shows that $|f_p(\alpha)|$ as a function of $\alpha$
is decreasing for $0<\alpha\le 1/2$ and increasing for
$1/2\le \alpha\le 1$. Thus
$$1-{2p\over (p^2+1)(p-1)}\le |f_p(\alpha)|\le 1,$$
where the lower bound holds true iff $\alpha=1/2$ and the upper bound
iff $\alpha=0$. It follows that $|A_{\chi}|\le 1$ with $A_{\chi}=1$ iff
$\chi$ is the principal character mod $d$.\\ 
\indent If $\chi'$ is the primitive Dirichlet character associated with $\chi$, 
then the Euler products of $A_{\chi}$ and $A_{\chi'}$ differ in at most
finitely many primes and hence can be simply related.\\
\indent It will be shown in Theorem \ref{intermsofachi} that $\rho(a,d)$ and $\delta(a,d)$ can
be expressed in terms of $A_{\chi}'s$, where $\chi$ ranges over the Dirichlet
characters mod $d$. The proof makes use of the following proposition.
\begin{Prop}
\label{puh}
Let $a\ge 1$. One has
$$\rho(a,d)={1\over \varphi(\delta)w\varphi(w)}
\prod_{p|\delta,~p\nmid w}\left(1-{1\over p(p-1)}\right)\prod_{p|\delta,~p|w}(1-{1\over p^2})$$
$$\sum_{\chi ({\rm mod~}\delta)}{\overline{\chi(\alpha)}}A_{\chi}
\prod_{p|w\atop p\nmid \delta}
{1+{\chi(p)-1\over p^2-\chi(p)}\over
1+{(\chi(p)-1)p\over (p^2-\chi(p))(p-1)}},$$
where $w=(a,d)$, $\alpha=a/w$ and $\delta=d/w$.
In particular, if $(a,d)=1$ then 
$$\rho(a,d)={1\over \varphi(d)}\prod_{p|d}
\left(1-{1\over p(p-1)}\right)\sum_{\chi ({\rm mod~}d)}{\overline{\chi(a)}}A_{\chi}.$$
\end{Prop}
\begin{Cor}
If $\gamma((a,d))(a,d)|d$, then
$$\rho(a,d)={\rho(\alpha,\delta)\over w\varphi(w)}
\prod_{p|(\delta,w)}{1-{1\over p^2}\over 1-{1\over p(p-1)}}.$$
\end{Cor}
{\it Proof of Proposition} \ref{puh}. From (\ref{nogonbedorven}) and Lemma \ref{nonmulttomult} one easily infers that
$$\rho(a,d)={1\over \varphi(\delta)w\varphi(w)}\sum_{\chi({\rm mod~}\delta)}{\overline{\chi(\alpha)}}
\sum_{v=1}^{\infty}{h_{\chi}(v)\varphi(w)\over v\varphi(vw)}.$$
On noting that the argument of the inner sum is multiplicative in $v$, 
the result follows on applying (\ref{muutje}) and Euler's product identity. \qed\\

\noindent {\tt Example 3}. 
One has $\rho(0,d)=1/(d\varphi(d))$ (in agreement with Proposition 
\ref{null}).
Let $\chi_3$ and $\chi_4$ denote the non-trivial character mod $3$, respectively mod $4$. One
finds $\rho(\pm 1,3)=5(1\pm A_{\chi_3})/12$ and $\rho(\pm 2,8)=5(1\pm A_{\chi_4})/12$.
Let $\chi$ be the character mod 5 
uniquely determined by $\chi(2)=i$.
One has $$\rho(3,5)={19\over 80}\left(1+2{\rm Re}(iA_{\chi})-A_{\chi^2}\right).$$ Using Table 2, these densities
can then be numerically approximated.\\

\noindent {\tt Example 4}. One has
$$\cases{\rho(2,8)=3\rho(1,4)/4 \cr \rho(6,8)=3\rho(3,4)/4}{\rm ~and~}
\cases{\rho(3,9)=8\rho(1,3)/45 \cr 
\rho(6,9)=8\rho(2,3)/45}{\rm ~and~}\cases{\rho(2,12)=3\rho(1,6)/4 \cr 
\rho(10,12)=3\rho(5,6)/4.}$$

\noindent One can now infer how the densities can be expressed in terms of $A_{\chi}$'s.
\begin{Thm}
\label{intermsofachi}
Let $\alpha$ and $\delta$ be as in Proposition {\rm \ref{puh}}.  Then
$$\rho(a,d)\in \mathbb Q(\zeta_{{\rm ord}_{\alpha}(\delta)})[A_{\chi_1},\cdots,A_{\chi_{\varphi(\delta)}}],$$
where $\chi_1,\cdots,\chi_{\varphi(\delta)}$ are the characters mod $\delta$.\\
\indent Let $\lambda$ denote 
Carmichael's function, that is $\lambda(d)$ equals the exponent of the group $(\mathbb Z/d\mathbb Z)^*$, then
$$\delta(a,d)\in \mathbb Q(\zeta_{\lambda(d)})[A_{\chi_1},\cdots,A_{\chi_{\varphi(d)}}].$$
\end{Thm}
{\it Proof}. The first part is a straightforward consequence of Proposition \ref{puh}. The
second part follows on applying part 7 of Theorem \ref{drie} and Proposition \ref{puh} together with the
observation that if $\delta|d$, then any character $\chi'$ mod $\delta$ can be lifted 
to a character
$\chi$ mod $d$, such that $A_{\chi'}=cA_{\chi}$, where $c\in \mathbb Q(\zeta_{\lambda(d)})$. \qed\\

\noindent The next result follows on combining Proposition \ref{puh} with
part 7 of Theorem \ref{drie}.
\begin{Prop}
\label{chiq}
Suppose that $q$ is a prime and $q\nmid a$. Then
$$\delta(a,q)={q^2-q-1\over (q-1)^2(q+1)}-{1\over (q-1)^2}
\sum_{\chi\ne \chi_0}\chi(-a)A_{\chi}.$$
\end{Prop}

\noindent The Euler product $A_{\chi}$ can also be expressed in terms
of $\rho(a,d)'s$. This yields
$$A_{\chi}={\sum_{a=1}^d \chi(a)\rho(a,d)\over 
\prod_{p|d}\left(1-{1\over p(p-1)}\right)}.$$
Thus, by (\ref{reducedtelop}) and $\rho(a,d)\ge 0$, one finds
$$|A_{\chi}|\le {\sum_{a=1,~(a,d)=1}^d \rho(a,d)\over 
\prod_{p|d}\left(1-{1\over p(p-1)}\right)}=1,$$
with equality iff $\chi$ is the principal character mod $d$.

\section{The numerical evaluation of $\delta(a,d)$ and $\rho(a,d)$}
Consider the numerical evaluation of the constant $A_{\chi}$. To this
end it turns out to be more convenient to consider
$$B_{\chi}:=\prod_{p}\left(1+{[\chi(p)-1]p\over [p^2-\chi(p)](p-1)}\right)=A_{\chi}\prod_{p|d}\left(1-{1\over p(p-1)}\right).$$ 
Recall that 
$L(s,\chi)$, the Dirichlet series for the character $\chi$, is defined,
for Re$(s)>1$ by $L(s,\chi)=\sum_{n=1}^{\infty}\chi(n)/n^s$.
\begin{Thm}
\label{numeriek}
Let $p_1(=2),p_2,\cdots$ be the sequence of consecutive primes. 
Let $\chi$ be any Dirichlet character and $n\ge 31$ $($hence $p_n\ge 127)$. Then
$$B_{\chi}=R_1AL(2,\chi)L(3,\chi)L(4,\chi)\prod_{k=1}^n 
\left(1+{\chi(p)\over p_k(p_k^2-p_k-1)}\right)(1-{\chi(p_k)\over p_k^3})(1-{\chi(p_k)\over p_k^4}),$$
with
$${1\over 1+p_{n+1}^{-3.85}}\le |R_1| \le 1+{1\over p_{n+1}^{3.85}}.$$
\end{Thm}
{\it Proof}. The first step is to note that 
$$B_{\chi}=AL(2,\chi)L(3,\chi)L(4,\chi)\prod_{k=1}^{\infty} 
\left(1+{\chi(p)\over p_k(p_k^2-p_k-1)}\right)(1-{\chi(p_k)\over p_k^3})(1-{\chi(p_k)\over p_k^4}),$$
An upper bound for the $k$th term in the latter product is given by
$$1+t^5{(2+2t+t^3+t^5)\over 1-t-t^2},$$ where $t=1/p_k$. For $t\ge 127$ some analysis shows that the
latter expression is bounded above by $1+t^{4.85}$. Using this one obtains
$$|R_1|\le \prod_{p>p_n}\left(1+{1\over p^{4.85}}\right)<1+\sum_{m>p_n}{1\over m^{4.85}}\le 
1+{1\over p_{n+1}^{4.85}}+\int_{p_{n+1}}^{\infty}{dt\over
t^{4.85}}\le 1+{1\over p_{n+!}^{3.85}}.$$ 
A similar argument allows one to obtain the lower bound. \qed\\

\noindent Since the Artin constant (see e.g. \cite{Niklasch})
and $L(2,\chi)$, $L(3,\chi)$ and
$L(4,\chi)$ can be each evaluated with high numerical accuracy, Theorem
\ref{numeriek} allows one to compute $A_{\chi}$ with high numerical
accuracy. Using Proposition \ref{puh} and part 7 of Theorem \ref{drie},
$\rho(a,d)$, respectively $\delta(a,d)$, can then be evaluated with
high numerical precision.\\
\indent A more straightforward, but
numerically much less powerful, approach in computing $\rho(a,d)$ and
$\delta(a,d)$, is to invoke part 7 of Theorem \ref{drie} and compute
$\rho(a,d)$ using the identity $\rho(a,d)=A\sum_{t\equiv a({\rm mod~}d)}r(t)$.
One has the following estimates.
\begin{Prop}
\label{laatste}
Let $x\ge 6$. One has
$$0<\rho(a,d)-A\sum_{t\equiv a({\rm mod~}d)\atop t\le x}r(t)<{1.28\over x}.$$
\end{Prop}
\begin{Cor}
\label{BLAH}
Let $a>0$ and $a+d\ge 6$, then
$$0<\rho(a,d)-Ar(a)<{1.28\over a+d}.$$
\end{Cor}
The most important ingredient of the proof will is the following lemma (the idea of which
was suggested to the author by Carl Pomerance \cite{P}).
\begin{Lem}
\label{carl}
For $x\ge 6$ one has
$$\sum_{n>x}{1\over n\varphi(n)}<{2.1\over x}.$$
\end{Lem}
{\it Proof}. Using that $\varphi(n)\ge \log(2n)/(n\log 2)$ and
that $\sum_{k\ge y}1/k^2<1.075/y$ for $y\ge 6$, one finds, for $x\ge 6$,
\begin{eqnarray*}
\sum_{t>x}{1\over t\varphi(t)}&=&\sum_{t>x}{1\over t^2}\sum_{d|t}
{|\mu(d)|\over \varphi(d)}\nonumber\\
&=&\sum_{d=1}^{\infty}{|\mu(d)|\over d^2\varphi(d)}\sum_{r> x/d}{1\over r^2}\nonumber\\
&\le &\zeta(2)\sum_{d>x/6}{|\mu(d)|\over d^2\varphi(d)}
+{1.075\over x}\sum_{d\le x/6}{|\mu(d)|\over d\varphi(d)}\nonumber\\
&\le &\zeta(2)\sum_{d>x/6}{\log 2d\over d^3\log 2}+
{1.075\over x}\sum_{d\le x/6}{|\mu(d)|\over d\varphi(d)}\nonumber\\
&\le & \zeta(2)\int_{[x/6]}^{\infty}{\log 2t\over t^3\log 2}dt
+{1.075\over x}{\zeta(2)\zeta(3)\over \zeta(6)}\nonumber\\
&\le & {\zeta(2)\over 4}{(2\log(2[x/6])+1)\over 
{[x/6]}^2 \log 2}+{1.075\over x}{\zeta(2)\zeta(3)\over \zeta(6)}.\nonumber\\
\end{eqnarray*}
For $x\ge 45000$ the latter upper bound is bounded above by $2.1/x$. 
After calculating $\sum_{n=1}^{\infty}1/(n\varphi(n))$ with
enough precision (see Remark 1) and
using
that $$\sum_{t>x}{1\over n\varphi(n)}=\sum_{n=1}^{\infty}{1\over n\varphi(n)}
-\sum_{t\le x}{1\over n\varphi(n)}
<2.20386-\sum_{t\le x}{1\over n\varphi(n)},$$ the result follows after 
verification in the range $(6,45000)$ (this verification is easily
seen to require only a finite amount of computation, cf. 
\cite[Lemma 4]{Moreetsjebbie}). \qed\\

\noindent The above argument can be easily adapted to show that
$$x\sum_{n>x}{1\over n\varphi(n)}\sim {\zeta(2)\zeta(3)\over \zeta(6)}
={1.943\cdots}.$$

\noindent {\it Proof of Proposition} \ref{laatste}. From (\ref{rhoaccent})
and $r(t)\ge 0$, one infers that
$$A\sum_{t\equiv a({\rm mod~}d)\atop t\le x}r(t)<\rho(a,d)\le \sum_{t\equiv a({\rm mod~}d)\atop t\le x}r(t)+
A\sum_{t>x}r(t)$$
By (\ref{er}) one has
$$A\sum_{t>x}r(t)\le {6\over \pi^2}\sum_{t>x}{1\over t\varphi(t)}.$$
On invoking Lemma \ref{carl}, the result then follows. \qed\\

\noindent The terms $AL(2,\chi)L(3,\chi)L(4,\chi)$
in Theorem \ref{numeriek} form the beginning of an expansion of $B_{\chi}$ in terms of special
values of L-series.
\begin{Thm}
\label{fiboconv}
Define numbers $G_{j+1}^{(r)}$ by
$${(-1)^r\over r}\sum_{d|r}{\mu(d)(-1)^{r\over d}\over (1-z^d-z^{2d})^{r/d}}=\sum_{j=0}^{\infty}G_{j+1}^{(r)}z^j.$$
One has
$$B_{\chi}=A{L(2,\chi)L(3,\chi)\over
L(6,\chi^2)}\prod_{r=1}^{\infty}\prod_{k=3r+1}^{\infty}L(k,\chi^r)^{\lambda(k,r)},$$
where $(-1)^{r-1}\lambda(k,r)=G_{k-3r+1}^{(r)}\in \mathbb Z_{>0}$.
\end{Thm}
Note that as formal series $(1-z-z^2)^{-1}=\sum_{j=0}^{\infty}F_{j+1}z^j$ where $F_j$ denotes the
$j$th Fibonacci number (thus $F_0=0$, $F_1=1$ etc.). The numbers defined by 
$(1-z-z^2)^{-r}=\sum_{j=0}^{\infty}F_{j+1}^{(r)}z^j$ are known as convolved Fibonacci numbers and
hence a reasonable term for the integers $G_{j+1}^{(r)}$ might be 
`convoluted convolved Fibonacci numbers'. For the convenience of 
the reader Table 1 gives a small sample of these numbers.\\
\indent The positivity of the numbers $(-1)^{r-1}\lambda(k,r)$ is established in
\cite{Pifi}, where the numbers $G_{j+1}^{(r)}$ are investigated. The argument uses 
Witt's dimension formula for free Lie algebras. The remaining part of Theorem \ref{fiboconv}
follows from the following more general result. 
\begin{Thm}
\label{algemener}
Suppose that $f(z)$ allows a formal power series in $z$ having only integer coefficients,
i.e. $f(z)=\sum_{j\ge 1}a(j)z^j$ with $a(j)\in \mathbb Z$. 
Let $g(z)=\sum_{j\ge 1}|a(j)|z^j$ and let $j_0\ge 0$ denote the smallest integer such that $a(j)\ne 0$.
Let
$$H^{(r)}(z)={1\over r}\sum_{d|r}\mu(d)f(z^d)^{r/d}=\sum_{j=0}^{\infty}h(j,r)z^j.$$
Then, as formal power series in $y$ and $z$, one has
\begin{equation}
\label{dubbelproduct}
1-yf(z)=\prod_{k=1}^{\infty}\prod_{j=kj_0}^{\infty}(1-z^jy^k)^{h(j,k)},
\end{equation}
Moreover, the numbers $h(j,k)$ are integers.\\
\indent Let $\epsilon>0$ be 
fixed. The identity {\rm (\ref{dubbelproduct})} holds for all complex numbers $y$ and $z$
with $g(|z|)y<1-\epsilon$ and $|z|<\rho_c$, where $\rho_c$ is the radius of convergence of the Taylor
series of $g$ around $z=0$. If, moreover, $\rho_c>1/2$, $g(1/2)<1$ and
$\sum_p g({1\over p})$ converges, then
\begin{equation}
\label{ggg}
\prod_p\left(1-\chi(p)f({1\over p})\right)= \prod_{k=1}^{\infty}\prod_{j=kj_0}^{\infty}L(j,\chi^k)^{-h(j,k)}.
\end{equation}
\end{Thm}
{\it Proof of Theorem \ref{fiboconv}}.  
Note that
$$(1-YX^2)\left({1+{(Y-1)X^2\over (1-YX^2)(1-X)}\over 1-{X^2\over 1-X}}\right)=
\left(1+{YX^3\over 1-X-X^2})\right).$$
By the first part of Theorem \ref{algemener} one infers that, as formal series,
$$(1+{YX^3\over 1-X-X^2})=
(1-YX^{3})^{-1}(1-Y^2X^6)\prod_{r=1}^{\infty}\prod_{k=3r+1}^{\infty}(1-X^kY^r)^{(-1)^rG_{k-3r+1}^{(r)}},$$
on noting that $G_1^{(r)}=1$ for $r=2$ and $G_1^{(r)}=0$ for $r\ge 2$.
Apply the second part with $f(z)=-z^3/(1-z-z^2)$ (and 
hence $g(z)=z^3/(1-z-z^2)$). The Taylor series for 
$g$ has radius of convergence $\rho_c=(\sqrt{5}-1)/2>0.5$. Note that $g(x)\le 1/2$ for all $0\le x\le 1/2$.
Furthermore, $\sum_p g(1/p)< \sum_p {4\over p^3}<\infty$. \qed\\

\noindent {\tt Remark 2}. The Dirichlet character is an example of a completely multiplicative function $h$, i.e. $h(nm)=h(n)h(m)$
for all natural numbers $n$ and $m$. If one defines $L(s,h)$ by $L(s,h)=\sum_{n=1}^{\infty}h(n)n^{-s}$, then under
the same conditions, one may replace $\chi$ in 
Theorem \ref{algemener} by any completely multiplicative function $h$ satisfying
$|h(n)|\le 1$.\\

\noindent {\tt Remark 3}. If $0<\rho_c\le 1/2$ and
$\sum_{p>1/\rho_c}g(1/p)$ converges, then an identity of the type (\ref{ggg}) still holds, but with
Dirichlet L-functions being replaced by partial Dirichlet L-functions. The idea is just to leave out
the local factor $1-\chi(p)f(1/p)$ for sufficiently many small primes $p$ and then proceed as before, cf.
the proof of Theorem 1 of \cite{Moreeconstant} (in the formulation of Theorem 1 there, replace
$p_{n_0}+1>1/\beta$ (a typo) by $p_{n_0+1}>\beta$).\\

\noindent {\tt Remark 4}. The conditions in the latter part of the theorem ensure that $f(z)=O(z^2)$ for
small $z$. This ensures on its turn that in the double product in (\ref{ggg}) only factors $L(j,\chi^k)$ 
with $k\ge 1$ and $j\ge 2k\ge 2$ appear.\\

\noindent The proof given here of Theorem \ref{fiboconv} rests on the following lemma. 
\begin{Lem}
\label{ontwikkeling}
Suppose that $f(X,Y)=\sum_{j,k}\alpha(j,k)X^jY^k$ with $\alpha(j,k)$ integers and $f(0,0)=0$. Then there are unique
integers $e(j,k)$ such that, as formal series, one has
$$1+f(X,Y)=\prod_{j=0}^{\infty}\prod_{k=0\atop (j,k)\ne (0,0)}^{\infty}(1-X^jY^k)^{e(j,k)}.$$
\end{Lem}
{\it Proof}. The term $X^{j_1}Y^{k_1}$ is said to be of lower weight than $X^{j_2}Y^{k_2}$ if $k_1<k_2$ or $k_1=k_2$ and
$j_1<j_2$. Suppose that $X^jY^k$ is the term of lowest weight appearing in $f(X,Y)$. Then consider
$(1+f(X,Y))(1-X^jY^k)^{-a(j,k)}$. This can be written as $1+g(X,Y)$ where all the coefficients of $g(X,Y)$ are
integers and the term of lowest weight in $g(X,Y)$ has strictly larger weight than the term of lowest weight in
$f(X,Y)$. Now iterate.\\
\indent It is not obvious from this argument that if one starts with a different weight ordering of the terms $X^jY^k$ we
end up with the same integers $e(j,k)$. Suppose that $h(X)$ has integer coefficients, then the coefficients $e(n)$
in $1+h(X)=\prod_{n=1}^{\infty}(1-X^n)^{e(n)}$ are unique, cf. 
\cite{Moreeconstant}. Hence, by setting $X=0$, respectively $Y=0$, one obtains
that $e(0,k)$, respectively $e(j,0)$ are uniquely determined. Setting $Y=X^m$ one obtains
that $1+f(X,X^m)=\prod_{n=1}^{\infty}(1-X^n)^{v(n)}$, where $v(n)$ is uniquely determined and
$v(2m)=e(2m,0)+e(m,1)+e(0,2)$. The uniqueness of $e(0,2)$, $e(2m,0)$ and $f(2m)$ then implies the uniqueness of
$e(m,1)$. The proof
will be completed by using induction. So suppose one has established that $e(j,k)$ with $k\le r$ for some $r\ge 1$ are uniquely determined.
Using that $v((r+2)m)=\sum_{k=0}^{r+2}e((r+2-k)m,k)$, one infers by
the induction hypothesis and using that $e(0,r+2)$ and $v((r+2)m)$ are uniquely 
determined, that $e(m,r+1)$ is uniquely determined. \qed\\

\noindent {\it Proof of Theorem \ref{algemener}}.  By M\"obius inversion and the definition of $H^{(r)}(z)$ one infers that
$$f(z)^r=\sum_{d|r}{r\over d}H^{({r\over d})}(z^d)=\sum_{d|r}{r\over d}h(j,{r\over d})
\sum_{j=0}^{\infty}z^{jd},$$
from which it is inferred that
$$\sum_{r=1}^{\infty}y^rf(z)^r=\sum_{k=1}^{\infty}\sum_{j=0}^{\infty}h(j,k)k\sum_{d=1}^{\infty}z^{jd}y^{kd}.$$
The latter identity with both sides divided out by $y$ can be rewritten as
$${f(z)\over 1-yf(z)}=\sum_{k=1}^{\infty}\sum_{j=0}^{\infty}{h(j,k)kz^jy^{k-1}\over 
1-z^jy^k}.$$
Formal integration of both sides with respect to $y$ gives
$$-\log(1-yf(z))=-\sum_{k=1}^{\infty}\sum_{j=0}^{\infty}h(j,k)\log(1-z^jy^k),$$
whence 
$$1-yf(z)=\prod_{k=1}^{\infty}\prod_{j=0}^{\infty}(1-z^jy^k)^{h(j,k)}.$$
On writing $f_1(z)=f(z)/z^{j_0}$ and $y_1=yz^{j_0}$ and expanding $1-y_1f_1(z)(=1-yf(z))$ in terms
of $y_1$ and $z$, it is then seen that (\ref{dubbelproduct}) holds.
The integrality of $h(j,k)$ follows by Lemma \ref{ontwikkeling}.\\
\indent The formal argument can be certainly made rigorous in the situation where
\begin{equation}
\label{convergentie}
\sum_{k=1}^{\infty}\sum_{j=j_0}^{\infty}|h(j,k)kz^jy^k|<\infty,
\end{equation}
where one is in the situation of absolute convergence and interchanges in order of summation are hence
allowed. Note that for $x\ge 0$, $g$ is an non-decreasing function of the real variable $x$. Now note that
$$\sum_{j=j_0}^{\infty}|h(j,r)z^j|\le {1\over r}\sum_{d|r}g(|z|^d)^{r/d}\le g(|z|)^r,$$
where the assumption that $|z|<\rho_c$ is being used. The double sum in (\ref{convergentie}) is thus majorized by
$\sum_{r=1}^{\infty}r g(|z|)^r|y|^r$ which in the given $(y,z)$ region converges.\\
\indent By a similar argument the convergence $\sum_p \sum_{r=1}^{\infty}r g(1/p)^r$, which is a 
consequence of the convergence of $\sum_p g(1/p)$ (one uses here that $g$ is non-decreasing as
a function of the real variable $x$ for $x\ge 0$ and that
$g(1/2)<1$), ensures the convergence of the triple product 
\begin{equation}
\label{tripleproduct}
\prod_p\left(1-\chi(p)f({1\over p})\right)=\prod_p \prod_{k=1}^{\infty}\prod_{j=kj_0}^{\infty}
\left(1-{\chi(p)^k\over p^j}\right)^{h(j,k)}.
\end{equation}
(From the theory of infinite products use that a product $\prod(1+\epsilon_v)$ is called absolutely convergent
if $\sum{\epsilon_v}$ is absolutely convergent and that in an absolutely convergent product the factors can be
reordered without changing its value.) On bringing the outer product over the primes $p$ to the inside and
using the Euler product for a Dirichlet L-series, the result then follows. \qed\\

\vfil\eject

\section{Tables}

\noindent {\it Explanation to Table} 1.
Table 1 gives some values of convoluted convolved Fibonacci numbers $G_j^{(r)}$. These
numbers are defined in Theorem \ref{fiboconv}.\\

\noindent {\it Explanation to Table} 2.
For every character $\chi$ of modulus $\le 12$, $A_{\chi}$ can be deduced from the table
below. In every case the value of $\chi$ is given (in at most
two arguments) such that $\chi$ is uniquely determined by this. If $\chi$ itself is not in the table, 
its complex conjugate $\bar \chi$ will be 
(in which case one has $A_{\chi}={\overline{A_{\bar \chi}}}$) 
or $\chi$ is the principal character (in which case $A_{\chi}=1$).
Although $A_{\chi}$ for $\chi$ not a primitive character can be easily
related to $A_{\chi'}$ with $\chi'$ a primitive character, for the
convenience of the reader the
numerical approximations to $A_{\chi}$ for the non-primitive
characters are listed as well.\\

\noindent {\it Explanation to Table} 3. An entry in a column having
as header the number $a$ and in a row starting with an integer $d$, 
respectively a $-$, gives the first
five decimal digits of $\delta(a,d)$, respectively $\delta(a+6,d)$.
If an entry is in a row labelled $\approx$, let $\delta(a,d)$ be the
entry directly above it. Then the number given equals 
$N_{-19}(a,d)(x)/\pi(x)$ with $x=2038074743$ (and hence $\pi(x)=10^8$).\\

\noindent  {\it Explanation to Table} 4. Similar to that of Table 2 (and with the same value of $x$).
In case $d=\infty$ one has $\delta(a,d)=Ar(a)$ for $a\ge 1$ and
$N'_g(a,\infty)(x)$ denotes the number of primes $p\le x$ with $v_p(g)=0$
such that $g$ has index equal to $a$. Here $x=1299709$ (and hence $\pi(1299709)=
10^5$).\\

\medskip
\medskip
\noindent {\bf Acknowledgement}. The author thanks Yves Gallot for
writing a $C^{++}$ program that was used to create Tables 3 and 4.
Furthermore, he thanks the referee for some helpful comments.

\vfill\eject

\centerline{{\bf Table 1:} Convoluted convolved Fibonacci numbers $G_j^{(r)}$}
\medskip
\begin{center}
\begin{tabular}{|c|c|c|c|c|c|c|c|c|c|c|c|}\hline
$r\backslash j$&1&2&3&4&5&6&7&8&9&10&11\\ \hline
\hline\hline
1&1&1&2&3&5&8&13&21&34&55&89\\ \hline
2&1&1&3&5&11&19&37&65&120&210&376\\ \hline
3&0&1&3&7&17&37&77&158&314&611&1174\\ \hline
4&0&1&3&10&25&64&146&331&710&1505&3091\\ \hline
5&0&1&4&13&38&102&259&626&1457&3287&7224\\ \hline
\end{tabular}  
\end{center}

\medskip
\medskip
\centerline{{\bf Table 2:} Numerical evaluation of $A_{\chi}$}
\medskip
\medskip
\begin{center}
\begin{tabular}{|c|c|c|c|}\hline
$d$&$\chi$&$\chi$&$A_{\chi}$\\ \hline
\hline\hline
3&$\chi(2)=-1~$&-&$+~0.173977122429634\cdots$\\ \hline
4&$\chi(3)=-1~$&-&$+~0.643650679662525\cdots$\\ \hline
5&$\chi(2)=~~i~$&-&$+~0.364689626478581\cdots$\\ \hline
-&-&-&$+i0.224041094424738\cdots$\\ \hline
5&$\chi(2)=-1~$&-&$+~0.129307938528080\cdots$\\ \hline
6&$\chi(5)=-1~$&-&$+~0.869885612148171\cdots$\\ \hline
7&$\chi(3)=e^{\pi i/3}$&-&$+~0.218769298429369\cdots$\\ \hline
-&-&-&$+i0.235418433356679\cdots$\\ \hline
7&$\chi(3)=e^{4\pi i/3}$&-&$+~0.212612780475062\cdots$\\ \hline
-&-&-&$-i0.145188986908610\cdots$\\ \hline
7&$\chi(3)=-1~$&-&$+~0.611324432919373\cdots$\\ \hline
8&$\chi(3)=~~1~$&$\chi(5)=-1$&$+~0.837998503129360\cdots$\\ \hline
8&$\chi(3)=-1~$&$\chi(5)=~~1$&$+~0.643650679662525\cdots$\\ \hline
8&$\chi(3)=-1~$&$\chi(5)=-1$&$+~0.603907856267167\cdots$\\ \hline
9&$\chi(2)=e^{\pi i/3}$&-&         $+~0.578815911632924\cdots$\\ \hline
-&-&-&$+i0.334468140016295\cdots$\\ \hline
9&$\chi(2)=e^{4\pi i/3}$&-&$+~0.250710892521489\cdots$\\ \hline
-&-&-&$-i0.207858981269346\cdots$\\ \hline
9&$\chi(2)=-1~$&-&$+~0.173977122429634\cdots$\\ \hline
10&$\chi(3)=~~i~$&-&$+~0.779414790379699\cdots$\\ \hline
-&-&-&$+i0.123970019579663\cdots$\\ \hline
10&$\chi(3)=-1~$&-&$+~         0.646539692640401\cdots$\\ \hline
11&$\chi(2)=e^{\pi i/5}$&-&$+~0.657644343795360\cdots$\\ \hline
-&-&-&                     $+i0.151998116640767\cdots$\\ \hline
11&$\chi(2)=e^{2\pi i/5}$&-&$+~0.373259555803500\cdots$\\ \hline
-&-&-&                     $+i0.208638808901506\cdots$\\ \hline
11&$\chi(2)=e^{3\pi i/5}$&-&$+~0.187051722258759\cdots$\\ \hline
-&-&-&                     $+i0.232381723173172\cdots$\\ \hline
11&$\chi(2)=-1~$&-&$+~0.184204262987186\cdots$\\ \hline
12&$\chi(5)=~~1~$&$\chi(7)=-1~$&$+~0.919500970946465\cdots$\\ \hline
12&$\chi(5)=-1~$&$\chi(7)=~~1~$&$+~0.869885612148171\cdots$\\ \hline
12&$\chi(5)=-1~$&$\chi(7)=-1~$&$+~0.841259078358102\cdots$\\ \hline
\end{tabular}
\end{center}

\vfill\eject
\centerline{{\bf Table 3:} $\delta(a,d)$ and approximation to $\delta_{-19}(a,d)$}
\medskip
\medskip
\begin{center}
\begin{tabular}{|c|c|c|c|c|c|c|}\hline
$a$&0&1&2&3&4&5\\ \hline
\hline\hline
$d=2$&0.66666&0.33333&-&-&-&-\\ \hline
$\approx$&0.66667&0.33333&-&-&-&-\\ \hline
3&0.37500&0.35599&0.26900&-&-&-\\ \hline
$\approx$&0.37502&0.35602&0.26897&-&-&-\\ \hline
4&0.33333&0.16666&0.33333&0.16666&-&-\\ \hline
$\approx$&0.33334&0.16664&0.33333&0.16669&-&-\\ \hline
5&0.20833&0.23542&0.17799&0.23400&0.14424&-\\ \hline
$\approx$&0.20831&0.23572&0.17829&0.23373&0.14395&-\\ \hline
$6$&0.25000&0.06067&0.12134&0.12500&0.29532&0.14766\\ \hline
$\approx$&0.25001&0.06067&0.12132&0.12501&0.29534&0.14765\\ \hline
7&0.14583&0.15968&0.15483&0.11905&0.16351&0.15567\\ \hline
$\approx$&0.14584&0.15965&0.15467&0.11915&0.16367&0.15573\\ \hline
-&0.10141&-&-&-&-&-\\ \hline
$\approx$&0.10129&-&-&-&-&-\\ \hline
8&0.16666&0.08333&0.16666&0.08333&0.16666&0.08333\\ \hline
$\approx$&0.16667&0.08332&0.16664&0.08335&0.16667&0.08332\\ \hline
-&0.16666&0.08333&-&-&-&-\\ \hline
$\approx$&0.16669&0.08334&-&-&-&-\\ \hline
9&0.12500&0.11866&0.08966&0.12500&0.11866&0.08966\\ \hline
$\approx$&0.12501&0.11866&0.08966&0.12501&0.11868&0.08965\\ \hline
-&0.12500&0.11866&0.08966&-&-&-\\ \hline
$\approx$&0.12500&0.11867&0.08965&-&-&-\\ \hline
10&0.13888&0.07196&0.14393&0.08172&0.06810&0.06944\\ \hline
$\approx$&0.13888&0.07197&0.14408&0.08159&0.06783&0.06944\\ \hline
-&0.16345&0.03405&0.15227&0.07613&-&-\\ \hline
$\approx$&0.16374&0.03421&0.15214&0.07612&-&-\\ \hline
11&0.09166&0.09890&0.09811&0.09904&0.09848&0.07170\\ \hline
$\approx$&0.09166&0.09889&0.09805&0.09904&0.09859&0.07180\\ \hline
-&0.09940&0.09303&0.09297&0.09523&0.06143&-\\ \hline
$\approx$&0.09939&0.09303&0.09297&0.09526&0.06133&-\\ \hline
12&0.12500&0.03033&0.06067&0.06250&0.14766&0.07383\\ \hline
$\approx$&0.12500&0.03033&0.06065&0.06251&0.14767&0.07382\\ \hline
-&0.12500&0.03033&0.06067&0.06250&0.14766&0.07383\\ \hline
$\approx$&0.12501&0.03035&0.06067&0.06249&0.14767&0.07383\\ \hline
\end{tabular}  
\end{center}
\vfill\eject

\vfill\eject
\centerline{{\bf Table 4:} $\rho(a,d)$ and approximation to $\rho_{65537}(a,d)$}
\medskip
\medskip
\begin{center}
\begin{tabular}{|c|c|c|c|c|c|c|}\hline
$a$&0&1&2&3&4&5\\ \hline
\hline\hline
$d=2$&0.50000&0.50000&-&-&-&-\\ \hline
$\approx$&0.49994&0.50006&-&-&-&-\\ \hline
3&0.16666&0.48915&0.34417&-&-&-\\ \hline
$\approx$&0.16662&0.48924&0.34414&-&-&-\\ \hline
4&0.12500&0.41091&0.37500&0.08908&-&-\\ \hline
$\approx$&0.12497&0.41097&0.37497&0.08909&-&-\\ \hline
5&0.05000&0.44143&0.31320&0.10036&0.09498&-\\ \hline
$\approx$&0.05000&0.44150&0.31322&0.10035&0.09494&-\\ \hline
6&0.08333&0.38955&0.31706&0.08333&0.09959&0.02710\\ \hline
$\approx$&0.08330&0.38966&0.31705&0.08331&0.09958&0.02709\\ \hline
7&0.02380&0.40253&0.29923&0.08966&0.08471&0.03881\\ \hline
$\approx$&0.02380&0.40263&0.29923&0.08962&0.08470&0.03881\\ \hline
-&0.06123&-&-&-&-&-\\ \hline
$\approx$&0.06122&-&-&-&-&-\\ \hline
8&0.03125&0.38569&0.30818&0.07380&0.09375&0.02521\\ \hline
$\approx$&0.03124&0.38577&0.30818&0.07380&0.09372&0.02521\\ \hline
-&0.06681&0.01528&-&-&-&-\\ \hline
$\approx$&0.06679&0.01529&-&-&-&-\\ \hline
9&0.01851&0.39347&0.29075&0.08696&0.07829&0.02983\\ \hline
$\approx$&0.01851&0.39356&0.29075&0.08694&0.07829&0.02981\\ \hline
-&0.06118&0.01738&0.02358&-&-&-\\ \hline
$\approx$&0.06117&0.01740&0.02358&-&-&-\\ \hline
10&0.02500&0.38063&0.30067&0.07141&0.08456&0.02500\\ \hline
$\approx$&0.02500&0.38071&0.30068&0.07140&0.08452&0.02500\\ \hline
-&0.06080&0.01253&0.02895&0.01041&-&-\\ \hline
$\approx$&0.06079&0.01254&0.02895&0.01041&-&-\\ \hline
11&0.00909&0.39040&0.28866&0.07722&0.07791&0.02698\\ \hline
$\approx$&0.00910&0.39047&0.28865&0.07721&0.07791&0.02698\\ \hline
-&0.05543&0.01768&0.02331&0.01418&0.01909&-\\ \hline
$\approx$&0.05541&0.01769&0.02331&0.01419&0.01908&-\\ \hline
12&0.02083&0.37819&0.29216&0.07231&0.07926&0.02170\\ \hline
$\approx$&0.02080&0.37827&0.29215&0.07230&0.07926&0.02169\\ \hline
-&0.06250&0.01136&0.02489&0.01101&0.02033&0.00540\\ \hline
$\approx$&0.06250&0.01139&0.02491&0.01101&0.02033&0.00539\\ \hline
$\cdots$&$\cdots$&$\cdots$&$\cdots$&$\cdots$&$\cdots$&$\cdots$\\ \hline
$\infty$&0.00000&0.37395&0.28046&0.06648&0.07011&0.01889\\ \hline
$\approx$&0.00000&0.37367&0.28124&0.06646&0.06913&0.01885\\ \hline
-&0.04986&0.00893&0.01752&0.00738&0.01417&0.00340\\ \hline
$\approx$&0.04962&0.00915&0.01796&0.00745&0.01449&0.00359\\ \hline
\end{tabular}  
\end{center}
\medskip
\vfill\eject

\medskip\noindent {\footnotesize KdV Institute,
Plantage Muidergracht 24, 1018 TV Amsterdam, The Netherlands.\\
e-mail: {\tt moree@science.uva.nl}}


\begin{thebibliography}{99}
\bibitem{Ballot} C. Ballot, Density of prime divisors of linear recurrences, {\it
 Mem. Amer. Math. Soc.} {\bf 115} (1995), no. 551, viii+102 pp..
\bibitem{CM4} K. Chinen and L. Murata, On a distribution property of
the residual order of $a({\rm mod~}p)$, math.NT/0211077.
\bibitem{Lenstra} {H. W. Lenstra, jr.}, {On Artin's conjecture and Euclid's
algorithm in global fields}, {\it Invent. Math.} {\bf 42} (1977), 201--224.
\bibitem{Early} P. Moree, On primes in arithmetic progression having a 
prescribed primitive root, 
{\it J. Number Theory} {\bf 78} (1999), 85--98.
\bibitem{Near} P. Moree, Asymptotically exact heuristics for (near) 
primitive roots. {\it J. Number Theory} {\bf 83} (2000), 155--181. 
\bibitem{Moreeconstant} { P. Moree}, Approximation of singular
series and automata, {\it Manuscripta Math.} {\bf 101} (2000), 385--399.
\bibitem{Moreetsjebbie} P. Moree, Chebyshev's bias for composite numbers
with restricted prime divisors, math.NT/0112100, to appear in Math. Comp. (electronically
already available from their site).
\bibitem{Moreealleen} P. Moree,
On the distribution of the order and index of
$g({\rm mod~}p)$ over residue classes, math.NT/0211259.
\bibitem{Pifi} P. Moree, Convoluted convolved Fibonacci numbers, preprint, 
vide http://staff.science.uva.nl/${\sim}$moree/preprints.html.
\bibitem{Moreealleen2} P. Moree,
On the distribution of the order and index of
$g({\rm mod~}p)$ over residue classes, II, in preparation.
\bibitem{Murato} L. Murata, A problem analogous to Artin's conjecture for primitive 
roots and its applications, {\it Arch. Math. (Basel)} {\bf 57} (1991), 555--565. 
\bibitem{CM5} L. Murata and K. Chinen, On a distribution property of
the residual order of $a({\rm mod~}p)$-II, math.NT/0211083.
\bibitem{Niklasch} G. Niklasch, http://www.gn-50uma.de/alula/essays/Moree/Moree.en.shtml
\bibitem{Odoni} R.W.K. Odoni,
A conjecture of Krishnamurthy on decimal periods and some allied problems,
{\it J. Number Theory} {\bf 13} (1981), 303--319.
\bibitem{Pappalardi} F. Pappalardi, On Hooley's theorem with weights, 
Number theory, II (Rome, 1995). {\it Rend. Sem. Mat. Univ. Politec. Torino}
{\bf 53} (1995), 375--388. 
\bibitem{P} C. Pomerance, e-mail, November 27th, 2002.
\bibitem{Prachar} K. Prachar, {\it Primzahlverteilung}, Springer, 
New York, 1957.
\bibitem{Wagstaff} {S.S. Wagstaff, jr.},
Pseudoprimes and a generalization of Artin's conjecture, {\it Acta Arith.}
{\bf 41} (1982), 141--150.
\bibitem{Wiertelak1} K. Wiertelak, On the density of some sets of primes.
IV, {\it Acta Arith.} {\bf 43} (1984), 177--190.
\bibitem{Wiertelak2} K. Wiertelak,
On the density of some sets of primes
$p$, for which $n|{\rm ord}_pa$, 
{\it Funct. Approx. Comment. Math.} {\bf 28}
(2000), 237--241.
\end{thebibliography}
\end{document}